\documentclass{llncs}
\usepackage[english]{babel}
\usepackage{amsmath}
\usepackage{amsfonts, mathrsfs}
\usepackage{amssymb}
\usepackage{verbatim}
\usepackage{graphicx}
\usepackage{color}
\usepackage[numbers]{natbib}

\newtheorem{defn}{Definition}

\renewcommand{\labelenumi}{\alph{enumi}.)}
\def\Indicator{\mathop{\hskip0pt{1}}\nolimits}

\title{PageRank in scale-free random graphs\thanks{This research is partially funded by the EU-FET Open grant NADINE (288956).}}

\author{Ningyuan Chen$^1$, Nelly Litvak$^2$, Mariana Olvera-Cravioto$^1$}

\institute{Columbia University, 500 W. 120th Street, 3rd floor, New York, NY 10027\and University of Twente, P.O.Box 217, 7500AE, Enschede, The Netherlands}

\begin{document}

\maketitle

\begin{abstract}
We analyze the distribution of PageRank on a directed configuration model and show that as the size of the graph grows to infinity it can be closely approximated by the PageRank of the root node of an appropriately constructed tree. This tree approximation is in turn related to the solution of a linear stochastic fixed point equation that has been thoroughly studied in the recent literature.

\end{abstract}

\section{Introduction}\label{sec:intro}

Google's PageRank proposed by Brin and Page~\cite{Brin98b} is arguably the most influential technique for computing centrality scores of nodes in a network. Numerous applications include graph clustering~\cite{Andersen06}, spam detection~\cite{Gyongyi04}, and citation analysis~\cite{Chen06citations,Waltman2010eigenfactor}. In this paper we analyze the power law behavior of PageRank scores in scale-free directed random graphs. 

In real-world networks, it is often found that the fraction of nodes with (in- or out-) degree $k$ is $\approx c_0k^{-\alpha-1}$, usually $\alpha\in(1,3)$, see e.g.~\cite{HofstadRG} for an excellent review of the mathematical properties of complex networks. More than ten years ago Pandurangan et al.~\cite{Pandurangan2006} discovered the interesting fact that PageRank scores also exhibit power laws, with the same exponent as the in-degree. This property holds for a broad class of real-life networks~\cite{Volkovich07}. 
In fact, the hypothesis that this always holds in power-law networks is plausible. 

However, analytical mathematical evidence supporting this hypothesis is surprisingly scarce. As one of the few examples, Avrachenkov and Lebedev~\cite{Avrachenkov06gn} obtained the power law behavior of average PageRank scores in a preferential attachment graph by using Polya's urn scheme and advanced numerical methods.

In a series of papers, Volkovich et al.~\cite{Litvak2007InternetMath,Volkovich07,Volkovich2010asymptotic} suggested an analytical explanation for the power law behavior of PageRank by comparing it to the endogenous solution of a stochastic fixed point equation (SFPE). The properties of this equation and the study of its multiple solutions has itself been an interesting topic in the recent literature \cite{Alsm_Mein_10b, Jelenkovic2010aap, Jel_Olv_12a, Jel_Olv_12b, Olvera_12a, Als_Dam_Men_12}, and is related to the broader study of weighted branching processes. The tail behavior of the endogenous solution, the one more closely related to PageRank, was given in \cite{Jelenkovic2010aap, Jel_Olv_12a, Jel_Olv_12b, Olvera_12a}, where it was shown to have a power law under many different sets of assumptions.  However, the SFPE does not fully explain the behavior of PageRank in networks since it implicitly assumes that the underlying graph is an infinite tree, an assumption that is not in general satisfied in real-world networks. 

\begin{figure}[h,t]
\centering
\includegraphics[scale = 0.5, bb = 70 430 570 630, clip]{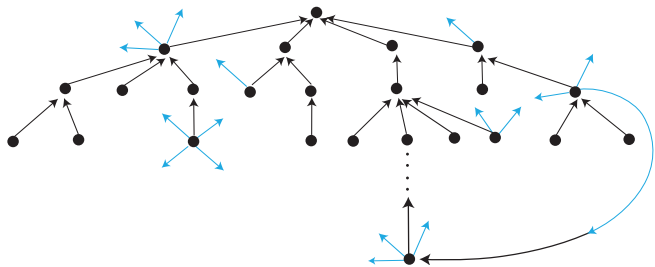}
\caption{Graph construction process. Unpaired outbound stubs are in blue.}\label{F.Tree}
\end{figure}

This paper makes a fundamental step further by extending the analysis of PageRank to graphs that are not necessarily trees. Specifically, we assume that the underlying graph is a directed configuration model (DCM) with given degree distributions, as developed by Chen and Olvera-Cravioto~\cite{chen2013directed}. We present numerical evidence that in this type of graphs the behavior of PageRank is very close to the one on trees. Intuitively, this is true for two main reasons: 1) the influence of remote nodes on the PageRank of an arbitrary node decreases exponentially fast with the graph distance; and 2) the DCM is asymptotically tree-like, that is, when we explore a graph starting from a given node, then with high probability the first loop is observed at a distance of order $\log n$, where $n$ is the size of the graph (see Figure \ref{F.Tree}). Our main result establishes analytically that PageRank in a DCM is well approximated by the PageRank of the root node of a suitably constructed tree as the graph size goes to infinity.

\vspace{-0.5cm}

\vspace{-0.3cm}

Section~\ref{sec:directed-random-graph} below describes the DCM as presented in~\cite{chen2013directed}. Then, in Section~\ref{sec:iterations} we analytically compare the PageRank scores in the DCM to their approximate value obtained after a finite number of power iterations. Next, in Section~\ref{sec:wbp} we explain how to couple the PageRank of a randomly chosen node with the root node of a suitable branching tree, and give our main analytical results. Finally, in Section~\ref{sec:numeric} we give numerical results validating our analytical work. The complete proofs for more general stochastic recursions, that also cover the PageRank case considered here, will be given in our upcoming paper~\cite{chen2014PR}.

%%% Directed Random Graphs
\section{Directed Random Graphs}\label{sec:directed-random-graph}

We will give below an algorithm, taken from \cite{chen2013directed}, that can be used to generate a scale-free directed graph.
Formally, power law distributions are modeled using the mathematical notion of regular variation. A nonnegative random variable $X$ is said to be regularly varying, if $\overline{F}(x):={P}(X>x)=L(x)x^{-\alpha}$, $x>0$, where  $L(\cdot)$ is a slowly varying function, that is, $\lim_{x\to\infty}L(tx)/L(x)=1$, for all $t>0$. 

Our goal now is to create a directed graph  $\mathcal{G}(n)$
with the property that the in-degrees and out-degrees will be approximately distributed, for large sizes of the graph, according to distributions $f^{\text{in}}_k  = P(\mathscr{N} = k)$, and $f^\text{out}_k = P(\mathscr{D} = k)$, $k = 0, 1, 2, 3,\dots$, respectively, where $E[\mathscr{N}] = E[\mathscr{D}]$. The only condition needed is that these distributions satisfy
$$\overline{F^{\text{in}}}(x) = \sum_{k > x} f_k^{\text{in}} \leq x^{-\alpha} L_{\text{in}}(x) \qquad \text{and} \qquad \overline{F^{\text{out}}} (x) =  \sum_{k > x} f_k^{\text{out}} \leq x^{-\beta} L_{\text{out}}(x), $$
for some slowly varying functions $L_\text{in}(\cdot)$ and $L_\text{out}(\cdot)$, and $\alpha, \beta > 1$.

The first step in our procedure is to generate an appropriate bi-degree sequence
\[
({\bf N}_n, {\bf D}_n) = \{ (N_i, D_i): 1 \leq i \leq n \}
\]
representing the $n$ nodes in the graph. The algorithm given below will ensure that the in- and out-degrees follow closely the desired distributions and also that the sums of in- and out-degrees are the same:
\[L_n:=\sum_{i=1}^nN_i=\sum_{i=1}^nD_i.\]
Denote
\[\kappa_0=\min\{1-\alpha^{-1}, 1-\beta^{-1}, 1/2\}.\]

{\bf Algorithm~1.} {\it Generation of a bi-degree sequence with given in-/out-degree distributions.}
\begin{enumerate} \itemsep 0pt
\renewcommand{\labelenumi}{\arabic{enumi}.}
\item Fix $0 < \delta_0 < \kappa_0$.
\item Sample an i.i.d. sequence $\{\mathscr{N}_1, \dots, \mathscr{N}_n\}$ from distribution $F^\text{in}$.

\item Sample an i.i.d. sequence $\{ \mathscr{D}_1, \dots, \mathscr{D}_n\}$ from distribution $F^\text{out}$, independent of $\{\mathscr{N}_i\}$.

\item Define $\Delta_n= \sum_{i=1}^n (\mathscr{N}_n- \mathscr{D}_n)$. If $|\Delta_n| \leq n^{1-\kappa_0 + \delta_0}$ proceed to step 5; otherwise repeat from step 2.

\item Choose randomly $|\Delta_n|$ nodes $\{i_1, i_2, \dots, i_{|\Delta_n|}\}$ without replacement and let
    \begin{align*}
        N_i &= \begin{cases}
        \mathscr{N}_i + 1 & \text{if $\Delta_n < 0$ and $i\in \{ i_1,i_2,\dots,i_{\Delta_n} \} $,}\\
        \mathscr{N}_i & \text{otherwise,}
    \end{cases}\\
    D_i &=\begin{cases}
        \mathscr{D}_i + 1 & \text{if $\Delta_n \ge 0$ and $i\in \{ i_1,i_2,\dots,i_{\Delta_n} \} $,}\\
        \mathscr{D}_i & \text{otherwise.}
    \end{cases}
  \end{align*}
\end{enumerate}

{\bf Remark:} It was shown in \cite{chen2013directed} that 
\begin{equation} \label{eq:Difference}
P\left( |\Delta_n| > n^{1-\kappa_0 + \delta_0} \right) = O\left(  n^{-\delta_0 (\kappa_0-\delta_0)/(1-\kappa_0) } \right)
\end{equation}
as $n \to \infty$, and therefore the Algorithm 1 will always terminate after a finite number of steps (i.e., it will eventually proceed to step 5). 

Having obtained a realization of the bi-degree sequence $({\bf N}_n, {\bf D}_n)$, we now use the configuration model to construct the random graph. The idea in the directed case is essentially the same as for undirected graphs. To each node $v_i$ we assign $N_i$ inbound half-edges and $D_i$ outbound half-edges; then, proceed to match inbound half-edges to outbound half-edges to form directed edges.
To be more precise, for each unpaired inbound half-edge of node $v_i$ choose randomly from all the available unpaired outbound half-edges, and if the selected outbound half-edge belongs to node, say, $v_j$, then add a directed edge from $v_j$ to $v_i$ to the graph; proceed in this way until all unpaired inbound half-edges are matched. Note that the resulting graph is not necessarily simple, i.e., it may contain self-loops and multiple edges in the same direction. 

We point out that conditional on the graph being simple, it is uniformly chosen among all simple directed graphs having bi-degree sequence $({\bf N}_n, {\bf D}_n)$ (see \cite{chen2013directed}). Moreover, it was also shown in \cite{chen2013directed} that, provided $\alpha, \beta > 2$, the probability of obtaining a simple graph through this procedure is bounded away from zero, and therefore one can obtain a simple graph having $({\bf N}_n, {\bf D}_n)$ as its bi-degree sequence by simply repeating the algorithm enough times. When we can only ensure that $\alpha, \beta > 1$, then a simple graph can still be obtained without loosing the distributional properties of the in- and out-degrees by erasing the self-loops and merging multiple edges in the same direction. These considerations about the graph being simple are nonetheless irrelevant to the current paper.

%%% Finitely many iterations
\section{PageRank iterations in the DCM}
\label{sec:iterations}

Although PageRank can be thought of as the solution to a system of linear equations, we will show in this section how it is sufficient to consider only a finite number of power iterations to obtain an accurate approximation for the PageRank of all the nodes in the graph. We first introduce some notation.

Let $M = M(n) \in \mathbb{R}^{n\times n}$ be matrix constructed as follows:
$$M_{i,j} = \begin{cases}
s_{ij} c/D_{i}, & \text{ if there are $s_{ij}$ edges from $i$ to $j$}, \\
0, & \text{ otherwise,}
\end{cases}$$
and let ${\bf 1}$ be the row vector of ones.  In the classical definition~\cite{LangvilleMeyer}, PageRank $\pi=(\pi_1,\ldots, \pi_n)$ is the unique solution to the following equation: 
\begin{equation}
\label{eq:pi}
\pi = \pi(c M) +\frac{1-c}{n}{\bf 1},
\end{equation}
where $c \in (0, 1)$ is a constant known as the damping factor. Rather than analyzing $\pi$ directly, we consider instead its scale-free version
\begin{equation}
\label{eq:PRdef}
n\pi= :{\bf R}=  {\bf R} (c M) +(1-c){\bf 1}
\end{equation}
obtained by multiplying \eqref{eq:pi} by the size of the graph $n$. Moreover, whereas $\pi_i$ is a probability distribution ($\pi_i \geq 0$ for all $i$ and $\pi {\bf 1}^T = 1$), its scale-free version ${\bf R} = (R_1, \dots, R_n)$ has components that are essentially unbounded for large $n$ and that satisfy $E[R_i]= 1$ for all $n$ and all $1\leq i \leq n$ (hence the name {\em scale-free}). 

One way to solve the system of linear equations given in \eqref{eq:PRdef} is via power iterations. We define the $k$th iteration of PageRank on the graph as follows. First initialize PageRank with a vector ${\bf r}_0 = r_0 {\bf 1}$, $r_0 \geq 0$, and then iterate according to ${\bf R}^{(n,0)} = {\bf r}_0$ and 
\begin{align*}
{\bf R}^{(n,k)}&= {\bf R}^{(n,k-1)} M + (1-c) {\bf 1} = (1-c){\bf 1} \sum_{i=0}^{k-1}  M^i + {\bf r}_0 M^k 
\end{align*}
for $k \geq 1$. In this notation, ${\bf R} = {\bf R}^{(n,\infty)}$, and our main interest is to analyze the distribution of the PageRank of a randomly chosen node in the DCM, say $R^{(n,\infty)}_1$. The first step of the analysis is to compare ${\bf R}^{(n,\infty)}$ to its $k$th iteration ${\bf R}^{(n,k)}$. To this end, note that ${\bf R}^{(n,\infty)} = (1-c){\bf 1} \sum_{i=0}^\infty  M^i$, and therefore,
$${\bf R}^{(n,k)} - {\bf R}^{(n,\infty)} = {\bf r}_0 M^k - (1-c){\bf 1} \sum_{i=k}^\infty  M^i .$$
Moreover,
\begin{align*}
\left|\left| {\bf R}^{(n,k)} - {\bf R}^{(n,\infty)} \right|\right|_1 &\leq \left|\left| {\bf r}_0 M^k \right|\right|_1 + (1-c)\sum_{i=0}^\infty \left|\left| {\bf 1} M^{k+i} \right|\right|_1 \\
&\leq r_0 n \left|\left| M^k \right|\right|_\infty + (1-c) n \sum_{i=0}^\infty \left|\left| M^{k+i} \right|\right|_\infty,
\end{align*}
where for the last inequality we used the observation that
$$\left|\left| {\bf 1} M^r \right|\right|_1 = \sum_{j=1}^n \sum_{i=1}^n (M^r)_{ij} = \sum_{i=1}^n \left|\left| (M^r)_{i\bullet} \right|\right|_1 \leq n \left|\left| M^r \right|\right|_\infty,$$
where $A_{i\bullet}$ denotes the $i$th row of matrix $A$.  Furthermore, since $M$ is equal to $c$ times an adjacency matrix, we have 
$$\left|\left| M^r \right|\right|_\infty \leq || M ||_\infty^r = c^r.$$
It follows that
\begin{align}
\label{eq:norm-k}
\left|\left| {\bf R}^{(n,k)} - {\bf R}^{(n,\infty)} \right|\right|_1 &\leq r_0 n c^k + (1-c) n \sum_{i=0}^\infty c^{k+i}= (r_0 +1) n c^k.
\end{align}

In general networks, the inequality for the $L_1$-norm (\ref{eq:norm-k}) does not provide information on convergence of specific coordinates and does not give a good upper bound for the quantity $|R_1^{(n,k)} - R_1^{(n,\infty)}|$ that we are interested in. However, the DCM has the additional property that all coordinates of the vector ${\bf R}^{(n,k)} - {\bf R}^{(n,\infty)}$ have the same distribution, since by construction, the DCM makes all permutations of the nodes' labels equally likely. This leads to the following observation. 

Let $\mathcal{F}_n = \sigma( ({\bf N}_n, {\bf D}_n))$ denote the sigma-algebra generated by the bi-degree sequence, which does not include information about the pairing process. Then, conditional on $\mathcal{F}_n$, 
\begin{align*}
E\left[ \left. \left| R^{(n,k)}_1 - R^{(n,\infty)}_1 \right| \right| \mathcal{F}_n \right] &= \frac{1}{n} E\left[ \left. \left|\left| {\bf R}^{(n,k)} - {\bf R}^{(n,\infty)} \right|\right|_1 \right| \mathcal{F}_n \right] \leq \left( r_0 + 1 \right) c^k,
\end{align*}
and Markov's inequality gives,
\begin{align}
P\left(  \left| R^{(n,\infty)}_1 - R^{(n,k)}_1\right| > \epsilon  \right) \notag &\leq E\left[ \epsilon^{-1} E\left[ \left. \left| R^{(n,k)}_1 - R^{(n,\infty)}_1 \right| \right| \mathcal{F}_n \right] \right] \notag \\
&\leq \left( r_0 + 1 \right) \epsilon^{-1} c^k \label{eq:PowerIterations}
\end{align}
for any $\epsilon > 0$.

Note that (\ref{eq:PowerIterations}) is a probabilistic statement, which is not completely analogous to (\ref{eq:norm-k}). In fact, (\ref{eq:PowerIterations}) states that we can achieve any level of precision with a pre-specified high probability by simply increasing the number of iterations $k$.  This leads to the following heuristic, that if the DCM looks locally like a tree for $k$ generations, where $k$ is the number of iterations needed to achieve the desired precision in \eqref{eq:PowerIterations}, then the PageRank of node 1 in the DCM will be essentially the same as the PageRank of the root node of a suitably constructed tree. The precise result and a sketch of the arguments will be given in the next section.

%%% Main Result: Coupling with a thorny branching tree
\section{Main Result: Coupling with a thorny branching tree}
\label{sec:wbp}

As mentioned in the previous section, we will now show how to identify $R^{(n,k)}_1$ with the PageRank of the root node of a tree. To start, we construct a variation of a branching tree where each node has an edge pointing to its parent but also has number of outbound stubs or half-edges that are pointing outside of the tree (i.e., to some auxiliary node).  We will refer to this tree as a Thorny Branching Tree (TBT), the name ``thorny" referring to the outbound stubs (see Figure \ref{F.Tree}). 

To construct simultaneously the graph $\mathcal{G}(n)$ and the TBT, denoted by $\mathcal{T}$, we start by choosing a node uniformly at random, and call it node 1 (the root node). This first node will have $N_1$ inbound stubs which we will proceed to match with randomly chosen outbound stubs. These outbound stubs are sampled independently and with replacement from all the possible $L_n=\sum_{i=1}^nD_i$ outbound stubs, discarding any outbound stub that has already been matched. This corresponds to drawing independently at random from the distribution 
\begin{align}
f_n(i,j) &= P(\text{node has $i$ offspring, $j$ outbound links }| \mathcal{F}_n ) \notag \\
&= \sum_{k=1}^n \Indicator(N_k = i, D_k = j)P(\text{an outbound stub of node $k$ is sampled }| \mathcal{F}_n) \notag \\
&=  \sum_{k=1}^n \Indicator(N_k = i, D_k = j) \frac{D_k}{L_n}. \label{eq:randomJointDistr}
\end{align}
This is a so-called size-biased distribution, since nodes with more outbound stubs are more likely to be chosen. 

To keep track of which outbound stubs have already been matched we will label them 1, 2, or 3 according to the following rule:
\begin{itemize}
\item [1.] Outbound stubs with label 1 are stubs belonging to a node that is not yet attached to the graph.
\item [2.] Outbound stubs with label 2 belong to nodes that are already part of the graph but that have not yet been paired with an inbound stub.
\item [3.] Outbound stubs with label 3 are those which have already been paired with an inbound stub and now form an edge in the graph.  
\end{itemize}

Let $Z_r$, $r\ge 0$, denote the number of inbound stubs of all the nodes in the graph at distance $r$ of the first node. Note that $Z_0=N_1$ and $Z_r$ is also the number of nodes at distance $(r+1)$ of the first node.

To draw the graph we initialize the process by labeling all outbound stubs with a 1, except for the $D_1$ outbound stubs of node 1 that receive a 2. We then start by pairing the first of the $N_1$ inbound stubs with a randomly chosen outbound stub, say belonging to node $j$. Then node $j$ is attached to the graph by forming an edge with node 1, and all the outbound stubs from the new node are now labeled 2. In case that $j=1$ the pairing forms a self-loop and no new nodes are added to the graph. Next, we label the chosen outbound stub with a 3, since it has already been paired, and in case $j \neq 1$, give all the other outbound stubs of node $j$ a label 2. We continue in this way until all $N_1$ inbound stubs of node 1 have been paired, after which we will be left with $Z_1$ unmatched inbound stubs that will determine the nodes at distance 2 from node 1.  In general, the $k$th iteration of this process is completed when all $Z_{k-1}$ inbound stubs have been matched with an outbound stub, and the process ends when all $L_n$ inbound stubs have been paired. Note that whenever an outbound stub with label 2 is chosen a cycle or double edge is formed in the graph. If at any point we sample an outbound stub with label 3 we simply discard it and do a redraw until we obtain an outbound stub with labels 1 or 2. 

We now explain the coupling with the TBT. We start with the root node (node 1, generation 0) that has $\hat N_1 = N_1$ offspring. Let $\hat Z_k$ denote the number of individuals in generation $k+1$ of the tree, $\hat Z_0 = \hat N_1$. For $k \geq 1$, each of the $\hat Z_{k-1}$ individuals in the $k$th generation will independently have offspring and outbound stubs according to the random joint distribution $f_n(i,j)$ given in \eqref{eq:randomJointDistr}. 

The coupling of the graph and the TBT is done according to the following rules: \vspace{-5pt}
\begin{itemize} 
\item [1.] If an outbound stub with label 1 is chosen, then both the graph and the TBT will connect the chosen outbound stub to the inbound stub being matched, resulting in a node being added to the graph and an offspring being born to its parent. In particular, if the chosen outbound stub corresponds to node $j$, then the new offspring in the TBT will have $D_j - 1$ outbound stubs (pointing to the auxiliary node) and $N_j$ inbound stubs (number of offspring).  We then update the labels by giving a 2 label to all the `sibling' outbound stubs of the chosen outbound stub, and a 3 label to the chosen outbound stub itself. 

\item [2.] If an outbound stub with label 2 is sampled it means that its corresponding node already belongs to the graph, and a cycle, self-loop, or multiple edge is created. In $\mathcal{T}$, we proceed as if the outbound stub had label~1 and create a new node, which is a copy of the drawn node. The coupling between DCM and TBT breaks at this point. 

\item [3.] If an outbound stub with label 3 is drawn it means that this stub has already been matched, and the coupling breaks as well. In $\mathcal{T}$, we again proceed as if the outbound stub had had a label~1. In the graph we do a redraw.
\end{itemize}

Note that the processes $Z_k$ and $\hat Z_k$ are identical as long as the coupling holds. Showing that the coupling holds for a sufficient number of generations is the essence of our main result.

\begin{defn}
Let $\tau$ be the number of generations in the TBT that can be completed before the first outbound stub with label 2 or 3 is drawn, i.e., $\tau = k$ iff the first inbound stub to draw an outbound stub with label 2 or 3 belonged to a node $i$, such that the graph distance between $i$ and the root node is exactly  $k$.
\end{defn}

The following result gives us an estimate as to when the coupling between the exploration process of the graph and the construction of the tree is expected to break. 

\begin{lemma} \label{L.CouplingBreaks}
Suppose $({\bf N}_n, {\bf D}_n)$ are constructed using Algorithm~1 with $\alpha > 1$, and $\beta>2$. Let $\mu = E[\mathscr{N}] = E[\mathscr{D}] > 1$. Then,  for any $1 \leq k \leq h \log n$ with $0 < h < 1/(2\log \mu)$ there exists a $\delta > 0$ such that,  
$$P\left( \tau \leq k  \right) = O\left( n^{-\delta} \right) \qquad \text{as } n \to \infty.$$
\end{lemma}

The proof of Lemma~\ref{L.CouplingBreaks} is rather technical, so we will only provide a sketch in this paper. The detailed proof will be given in \cite{chen2014PR}. 

\begin{proof}[Qualitative argument]
Let $\hat V_s$ be the number of outbound stubs of all nodes in generation $s$ of the tree. The intuition behind the proof is that for all $s = 1, 2, \dots$, neither $\hat Z_s$, nor $\hat V_s$ are expected to be much larger than their means:
$$E\left[ \left. \hat Z_s \right| \mathcal{F}_n\right] \approx \mu^{s+1} \qquad \text{and} \qquad E\left[ \left. \hat V_s \right| \mathcal{F}_n\right] \approx \lambda \mu^{s},$$
where $\lambda = E[\mathscr{D}^2]/\mu$. 
Next, note that an inbound stub of a node in the $r$th generation will be the first one to be paired with an outbound stub having label 2 or 3 with a probability bounded from above by
\[ 
P_r := \frac{1}{L_n} \sum_{s=0}^r \hat V_s\approx \frac{\lambda\mu^{r}}{n(\mu-1)}.\] 
Furthermore, for event $\{\tau = r\}$  to occur one of the $\hat Z_r$ inbound stubs must have been paired with an outbound stub with labels 2 or 3, which is bounded by the probability that a Binomial random variable with parameters $(\hat Z_r, P_r)$ is greater or equal than 1. Since $P_r=o(1)$ for $r\le k$, this probability is 
$1-(1-P_r)^{\hat Z_r}= P_r{\hat Z_r}(1+O(P_r)) = O\left( \mu^{2r}n^{-1} \right)$.

Formally, to ensure that the approximations given above are valid, we first show that the event
$$E_k = \left\{ \max_{0 \leq r \leq k} \mu^{-r}{\hat Z_r}{\mu^r} \leq x_n, \, \max_{0 \leq r \leq k} {\mu^{-(r+1)}} \sum_{s=0}^r \hat V_s \leq x_n \right\}$$
occurs with high probability as $n \to \infty$ for a suitably chosen $x_n \to \infty$. Then, sum over $r = 0, 1, \dots, k$ the events $\{ \tau =r, \, E_k\}$ to obtain that $P(\tau \leq k, \, E_k) = O\left( \mu^{2k} n^{-1}\right)$, which goes to zero for $k \leq h \log n$. 
\end{proof}

Our main result is now a direct consequence of the bound derived in \eqref{eq:PowerIterations} and Lemma \ref{L.CouplingBreaks} above, since before the coupling breaks $R_1^{(n,k)}$ and the PageRank, computed after $k$ iterations, of the root node of the coupled tree coincide.

\begin{theorem} \label{T.Main}
Suppose $({\bf N}_n, {\bf D}_n)$ are constructed using Algorithm~1 with $\alpha > 1$, and $\beta>2$. Let $\mu = E[\mathscr{N}] = E[\mathscr{D}] > 1$ and $c \in (0,1)$. Then, for any $\epsilon > 0$ and any $1 \leq k \leq h \log n$ with $0 < h < 1/(2\log \mu)$  there exists a $\delta > 0$ such that,
$$P\left( \left| R_1^{(n,\infty)} - \hat R_1^{(n,k)} \right| > \epsilon \right) \leq (r_0+1) \epsilon^{-1} c^k + O\left( n^{-\delta}  \right),$$
as $n \to \infty$, where $\hat R_1^{(n,k)}$ is the PageRank, after $k$ iterations, of the root node of the TBT described above. 
\end{theorem}

In the forthcoming paper \cite{chen2014PR} we explore further the distribution of the PageRank of the root node of $\mathcal{T}$ and show that $\hat R_1^{(n,k)}$ converges to the endogenous solution of a SFPE on a weighted branching tree, as originally suggested in \cite{Litvak2007InternetMath,Volkovich07,Volkovich2010asymptotic}. Moreover, the tail behavior of this solution has been fully described in \cite{Volkovich2010asymptotic, Jelenkovic2010aap, Jel_Olv_12a}.

%%% Numerical Results
\section{Numerical Results}
\label{sec:numeric}

In this last section we give some numerical results showing the accuracy of the TBT approximation to the PageRank in the DCM. To generate the bi-degree sequence we use as target distributions two Pareto-like distributions. More precisely, we set
\begin{equation*}
    \mathscr N_i=\lfloor X_{1,i}+Y_{1,i} \rfloor, \quad \mathscr D_i=\lfloor X_{2,i}+Y_{2,i} \rfloor,
\end{equation*}
where the $\{X_{1,i}\}$ and the $\{X_{2,i}\}$ are independent sequences of i.i.d. Pareto random variables with shape parameters $\alpha > 1$ and $\beta > 2$, respectively, and scale parameters $x_1 = (\alpha-1)/\alpha$ and $x_2 = (\beta-1)/\beta$, respectively (note that $E[X_{1,i}] = E[X_{2,i}] = 1$ for all $i$). The sequences $\{Y_{1,i}\}$ and $\{Y_{2,i}\}$ are independent sequences, each consisting of i.i.d. exponential random variables with means $1/\lambda_1 > 0$ and $1/\lambda_2$, respectively. The addition of the exponential random variables allows more flexibility in the modeling of the in- and out-degree distributions while preserving a power law tail behavior; the parameters $\lambda_1, \lambda_2$ are also used to match the means $E[\mathscr{N}]$ and $E[\mathscr{D}]$. 

Once the sequences $\{\mathscr{N}_i\}$ and $\{\mathscr{D}_i\}$ are generated, we use Algorithm 1 to obtain a valid bi-degree sequence $({\bf N}_n, {\bf D}_n)$.  Given this bi-degree sequence we next proceed to construct the graph and the TBT simultaneously, according to the rules described in Section \ref{sec:wbp}. To compute ${\bf R}^{(n,\infty)}$ we perform power iterations with $r_0 = 1$ until $\| {\bf R}^{(n,k)}- {\bf R}^{(n,k-1)}\|_2<\epsilon_0$ for some tolerance $\epsilon_0$. We only generate the TBT for the required number of generations in each of the examples; the computation of $\hat R_1^{(n,k)}$ can be done recursively starting from the leaves using 
\begin{equation}\label{eq:pagerank-wbp}
    \hat R_i^{(n,0)}=1,\quad \hat R_i^{(n,k)}=\sum_{j\rightarrow i}c\hat R_j^{(n,k-1)}+(1-c), \quad k > 0,
\end{equation}
where $j \to i$ means that node $j$ is an offspring of node $i$.

Tables \ref{tab:numerical-1}-\ref{tab:numerical-3} below compare the PageRank of node 1 in the graph, $R_1^{(n,\infty)}$, the PageRank of node 1 only after $k$ power iterations, $R_1^{(n,k)}$, and the PageRank of the root node of the coupled tree after the same $k$ generations, $\hat R_1^{(n,k)}$. The magnitude of the mean squared errors (MSEs), computed using $R_1^{(n,\infty)}$ as the true value, is also given in each table. The tolerance for computing $R_1^{(n,\infty)}$ is set to $\varepsilon_0=10^{-6}$. For each $n$, we generate $100$ realizations of $\mathcal{G}(n)$ as well as of the corresponding TBTs and take the empirical average of the PageRank values and of the MSEs. Table \ref{tab:numerical-1} includes results for different sizes of the graph, and uses $k_n = \lfloor \log n\rfloor$ iterations for the finite approximations. We note that all the MSEs clearly decrease as $n$ increases since $k_n$ also increases with $n$. 

\vspace{-0.3cm}
\begin{table}[h]
    \centering
    \begin{tabular}{|c|c|c|c|c|c|}
		\hline
        $n$ & \hspace{5pt} $R_1^{(n,\infty)}$ & \hspace{5pt} $R_1^{(n,k_n)}$ & \hspace{5pt} $\hat R_1^{(n,k_n)}$ & MSE for $R_1^{(n,k_n)}$ & MSE for $\hat R_1^{(n,k_n)}$ \\
        \hline
        10 & 0.931 & 0.946 & 0.983 & 3.90E-03 & 4.20E-02 \\
        100 & 1.023 & 1.027 & 1.068 & 1.80E-04 & 3.70E-02 \\
        1000 & 1.000 & 1.002 & 1.010 & 1.20E-05 & 8.00E-04 \\
        10000 & 0.964 & 0.965 & 0.962 & 1.00E-06 & 7.50E-04 \\
				\hline
    \end{tabular}
		\medskip
    \caption{$\alpha=2$, $\beta=2.5$, $\lambda_1=1$, $c=0.5$, $k_n=\lfloor \log n \rfloor$.}
    \label{tab:numerical-1}
\end{table}

\vspace{-0.7cm}

Table~\ref{tab:numerical-2} illustrates the impact of using different values of $k$, with the error between $R_1^{(n,k)}$ and $R_1^{(n,\infty)}$ clearly decreasing as $k$ increases. The simulations were run on a graph with $n = 10,000$ nodes. We also point out that although the accuracy of finitely many PageRank iterations improves as $k$ gets larger, the MSE of the tree approximation seems to plateau after a certain point. In order to obtain a higher level of precision we also need to increase the size of the graph (as suggested by Theorem \ref{T.Main}). 

\vspace{-0.5cm}

\begin{table}[h]
    \centering
    \begin{tabular}{|c|c|c|c|c|c|}
		\hline
       \hspace{2pt} $k_n$ \hspace{3pt} & \hspace{5pt} $R_1^{(n,\infty)}$ & \hspace{5pt} $R_1^{(n,k_n)}$ & \hspace{5pt} $\hat R_1^{(n,k_n)}$ & MSE for $R_1^{(n,k_n)}$ & MSE for $\hat R_1^{(n,k_n)}$ \tabularnewline
        \hline
        2 & 0.908 & 0.933 & 0.928 & 7.1E-03 & 8.59E-03\tabularnewline
        4 & 0.929 & 0.933 & 0.933 & 1.5E-04 & 2.20E-04\tabularnewline
        6 & 0.908 & 0.909 & 0.910 & 5.4E-06 & 5.08E-05\tabularnewline
        8 & 0.883 & 0.884 & 0.884 & 8.8E-08 & 1.20E-06\tabularnewline
        10 & 0.948 & 0.949 & 0.950 & 7.6E-09 & 8.16E-05\tabularnewline
        15 & 0.932 & 0.932 & 0.932 & 7.9E-13 & 2.89E-05\tabularnewline\hline
    \end{tabular}
		\medskip 
		
    \caption{$n=10000$, $\alpha=2$, $\beta=2.5$, $\lambda_1=1$, $c=0.5$.}
    \label{tab:numerical-2}
\end{table}

\newpage

Table~\ref{tab:numerical-3} shows the same comparison as in Table~\ref{tab:numerical-2}, for fixed $n$, for different values of the damping factor $c$. As $c$ gets larger, the approximations provided by both $R_1^{(n,k_n)}$ and $\hat R_1^{(n,k_n)}$ get worse due to the slower convergence of PageRank.

\vspace{-0.3cm}
\begin{table}[h]
    \centering
    \begin{tabular}{|c|c|c|c|c|c|}
		\hline
               \hspace{2pt} $c$ \hspace{3pt} & \hspace{5pt} $R_1^{(n,\infty)}$ & \hspace{5pt} $R_1^{(n,k_n)}$ & \hspace{5pt} $\hat R_1^{(n,k_n)}$ & MSE for $R_1^{(n,k_n)}$ & MSE for $\hat R_1^{(n,k_n)}$ \tabularnewline        \hline
        0.1 & 1.011 & 1.011 & 1.011 & 3.8E-22 & 3.33E-09\tabularnewline
        0.3 & 0.958 & 0.958 & 0.958 & 9.8E-13 & 1.91E-07\tabularnewline
        0.5 & 0.898 & 0.898 & 0.899 & 2.7E-08 & 2.63E-06\tabularnewline
        0.7 & 0.755 & 0.757 & 0.760 & 2.4E-05 & 2.03E-04\tabularnewline
        0.9 & 0.663 & 0.764 & 0.799 & 8.3E-02 & 1.25E-01\tabularnewline\hline
    \end{tabular}
		\medskip
		
    \caption{$n=10000$, $\alpha=2$, $\beta=2.5$, $\lambda_1=1$, $k_n = \lfloor \log n \rfloor =9$.}
    \label{tab:numerical-3}
\end{table}

%The last table, Table~\ref{tab:numerical-4} gives the corresponding comparison when we vary the expected number of degrees, $\mu = E[\mathscr{N}] = E[\mathscr{D}]$. In this case the pattern of the MSEs is not so clear, but we would expect that the larger $\mu$ is the faster the coupling with the tree will break, and this should have an impact on the quality of the tree approximation for a fixed number of nodes $n$.  

%\begin{table}
%    \centering
%    \begin{tabular}[h]{|c|c|c|c|c|c|}
%		\hline
%               \hspace{2pt} $\mu$ \hspace{3pt} & \hspace{5pt} $R_1^{(n,\infty)}$ & \hspace{5pt} $R_1^{(n,k_n)}$ & \hspace{5pt} $\hat R_1^{(n,k_n)}$ & MSE for $R_1^{(n,k_n)}$ & MSE for $\hat R_1^{(n,k_n)}$ \tabularnewline     
%        \hline
%        1.5 & 0.975 & 0.975 & 0.975 & 2.3E-07 & 1.54E-05\tabularnewline
%        2.0 & 0.923 & 0.923 & 0.924 & 2.1E-08 & 1.21E-05\tabularnewline
%        2.5 & 0.945 & 0.945 & 0.948 & 2.1E-08 & 1.33E-04\tabularnewline
%        3.0 & 0.968 & 0.968 & 0.972 & 1.2E-08 & 1.01E-04\tabularnewline
%        3.5 & 0.892 & 0.892 & 0.897 & 5.4E-09 & 2.01E-04\\\hline
%    \end{tabular}
%				\medskip
%				
%    \caption{$n=10000$, $\alpha=2$, $\beta=2.5$, $c=0.5$ and $k_n= \lfloor \log n \rfloor = 9$.}
%    \label{tab:numerical-4}
%\end{table}

\vspace{-0.5cm}

\begin{figure}[!hb]
   \centering
    \includegraphics[scale=0.5, bb = 0 0 370 250, clip]{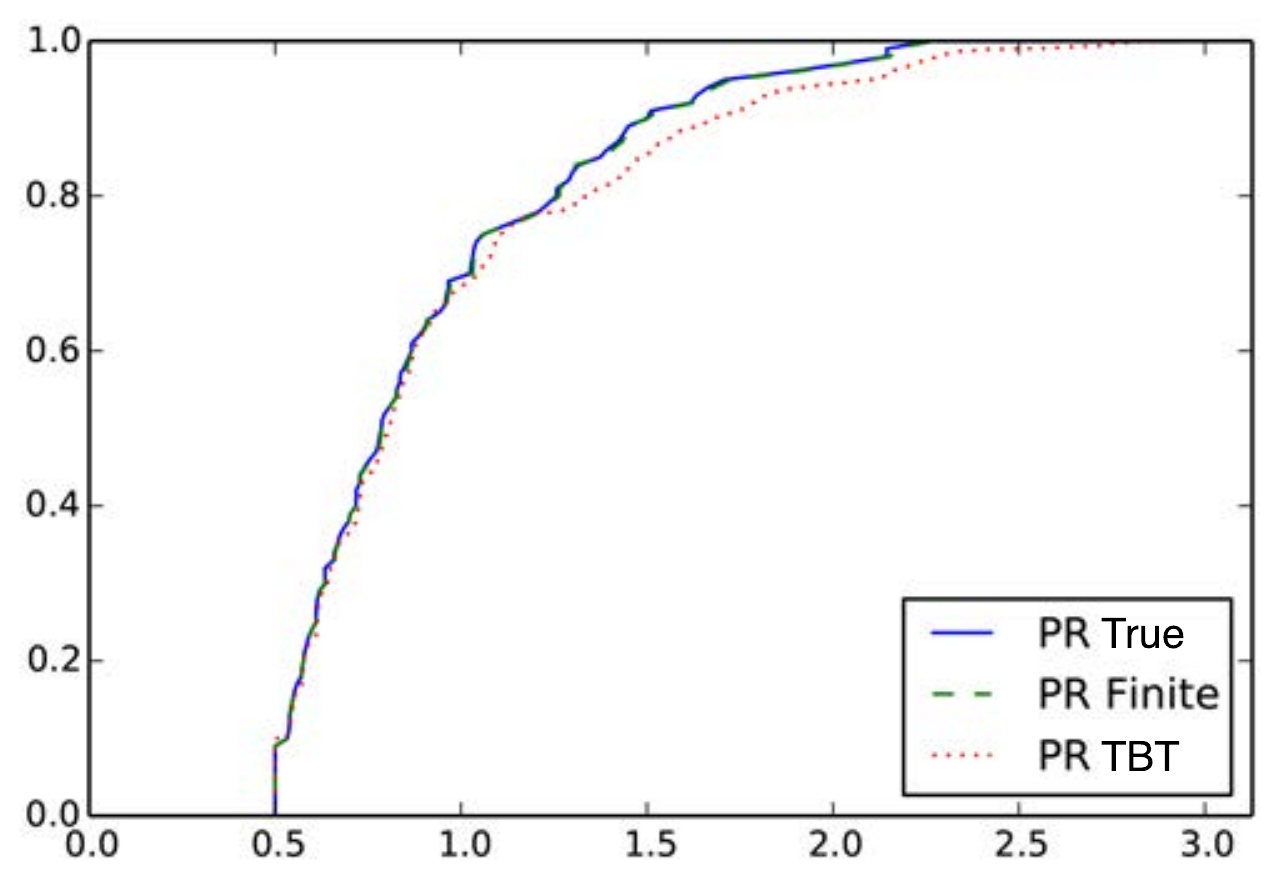}
    \caption{The empirical distributions of PageRank on $\mathcal{G}(n)$ (true and after finitely many iterations) and the empirical distribution of the PageRank of the root in the TBT.}
    \label{fig:empirical-pr}
\end{figure}

Our last numerical result shows how the distribution of PageRank on the TBT approximates the distribution of PageRank on the DCM. To illustrate this we generated a graph with $n = 100$ nodes and parameters $\alpha = 2$, $\beta = 2.5$, $\mu = 3$ and $c = 0.5$. We set the number of PageRank iterations (number of generations in the TBT) to be $k = 4$. We then computed the empirical CDFs of the PageRank of all nodes in the graph and that of the PageRank after only $k$ iterations. We also generated the coupled TBT 1000 times based on the same graph; each time by randomly choosing  some  node $i$ to be the root and computing $\hat R_i^{(n,k)}$ according to \eqref{eq:pagerank-wbp}. Figure \ref{fig:empirical-pr} plots  the empirical CDF of PagerRank on $\mathcal{G}(n)$, the empirical CDF of PageRank on $\mathcal{G}(n)$ after only $k$ iterations, and the empirical CDF of the PageRank of the 1000 root nodes after the same $k$ iterations. We can see that the CDFs of PageRank on $\mathcal{G}(n)$ after a finite number of iterations and that of the true PageRank on $\mathcal{G}(n)$ are almost indistinguishable. The PageRank on the TBT also approximates this distribution quite well, especially considering that $n = 100$ is not particularly large.

\bibliographystyle{plain}
\bibliography{ref}

\end{document}